\documentclass[11pt]{article}
\usepackage[margin=1in]{geometry}
\usepackage[T1]{fontenc}
\usepackage{lmodern}
\usepackage{microtype}
\usepackage{amsmath, amssymb, amsthm}
\usepackage{booktabs}
\usepackage{comment}
\usepackage{graphicx}
\usepackage{xcolor}
\usepackage{tikz}
\usetikzlibrary{positioning, arrows.meta, calc}
\usepackage{algorithm}
\usepackage{algpseudocode}
\usepackage[round, authoryear]{natbib}
\usepackage[hypertexnames=false, colorlinks, linkcolor=blue!50!black, citecolor=blue!50!black, urlcolor=blue!50!black]{hyperref}

\newtheorem{proposition}{Proposition}
\newtheorem{lemma}[proposition]{Lemma}
\newtheorem{theorem}[proposition]{Theorem}
\newcommand{\cmax}{C_{\max}}
\newcommand{\cpw}{\textsc{CPW}}
\newcommand{\est}{\mathrm{est}}
\newcommand{\lct}{\mathrm{lct}}
\newcommand{\ect}{\mathrm{ECT}}

\title{Fragment-Free Constraint Programming for the Preemptive Flexible Job Shop Scheduling Problem}
\author{Vincent Boyer\thanks{Tecnol\'ogico de Monterrey, School of Engineering and Sciences, Campus Monterrey, Mexico. E-mail:
\href{mailto:vincent.boyer@tec.mx}{vincent.boyer@tec.mx}.}
\and
M. Ang\'elica Salazar-Aguilar\thanks{ Universidad Aut\'onoma de Nuevo Le\'on (UANL), Facultad de Ingenier\'{\i} a Mec\'anica y El\'ectrica (FIME), Mexico. E-mail:
\href{mailto:maria.salazaragl@uanl.edu.mx}{maria.salazaragl@uanl.edu.mx}. Corresponding author.}}
\date{November 2026}

\begin{document}
\maketitle

\begin{abstract}
\noindent
This article addresses the preemptive flexible job shop scheduling problem with makespan minimization ($FJm \mid pmtn \mid C_{\max}$). We propose an exact constraint programming approach that completely avoids representing the processing fragments of an operation. Operations are placed in flexible time windows, where Horn's workload condition guarantees feasible preemptive scheduling, and Jackson's preemptive earliest-deadline-first rule reconstructs exact processing segments post-optimization in $O(n \log n)$ time per machine. To manage undecided machine assignments, we extend the \texttt{PreemptiveNoOverlap} global constraint to optional intervals within IBM ILOG CP Optimizer, performing overload checking and preemptive edge-finding in $O(n \log n)$ time while eliminating infeasible machine alternatives. The CP model is warm-started using a multi-round MILP assignment relaxation with Carlier's head--tail bounds and a window-based greedy heuristic. Extensive computational experiments on 276 pFJSSP benchmark instances and 78 pJSSP instances demonstrate that the proposed method significantly advances the state of the art. Within a 600~s time limit, it improves the upper bounds of the previously published logic-based Benders decomposition by 3.6\% on average, increases the number of proven optimal pFJSSP instances from 102 to 180, and closes 74 instances for the first time (including the entire Barnes family). On the pJSSP under identical hardware conditions, it proves 67 of 78 instances optimal (compared to 51 by the baseline solver Mistral) and closes 9 of the 11 long-standing open problems.

\end{abstract}

\section{Introduction}
\label{sec:intro}

The Job Shop Scheduling Problem (JSSP) is a classical combinatorial optimization problem involving a set of jobs and a set of machines. Each job consists of a sequence of operations with predetermined machine assignments and processing times. Operations must be processed without interruption, and each machine can handle at most one operation at a time. The objective is to determine the ordering and timing of operations that satisfy all precedence and machine-capacity constraints while minimizing the makespan, i.e., the completion time of the last job. In the flexible job shop scheduling problem (FJSSP), introduced by \citet{Brucker1990}, each operation may be processed by any machine of a set of eligible machines, with machine-dependent processing times, so that the problem combines an assignment and a sequencing decision. The FJSSP models many manufacturing systems with multi-purpose machines and has generated a vast literature, a recent survey can be found in \citet{DauzerePeres2024}.

In most of the existing literature, operations are not interruptible. Many settings, however, allow \emph{preemption}: an operation may be stopped and resumed later on the same machine, without penalty. This is the case when the resources are computers, when an urgent operation may interrupt a long one, or when the processing of a part can be split across shifts. Preemption can only shorten the optimal makespan, and the gain may be substantial \citep{Creemers2019}. It also changes the nature of the problem: the classical disjunctive reasoning of the job shop no longer applies, and the preemptive job shop is NP-hard already with two machines and three jobs, whereas its non-preemptive counterpart is then polynomially solvable \citep{Brucker1999}.

This paper addresses the \emph{preemptive flexible job shop scheduling problem} (pFJSSP) with makespan minimization, denoted $FJm\,|\,pmtn\,|\,\cmax$. Exact methods for this problem are recent. \citet{Juvin2023COR} proposed a time-indexed mixed-integer program, a constraint programming (CP) model in which each operation is split into unit-duration parts, and a logic-based Benders decomposition (LBBD) whose master problem assigns the operations to the machines and whose subproblem is a preemptive job shop solved by branch and bound. They introduced a benchmark of 276 instances, on which the LBBD clearly dominates the other approaches, but proves the optimality of only 102 instances within one hour. For the pJSSP, \citet{Juvin2023CP} showed that the fragments of the operations need not be represented at all: if each operation is given a time window, a preemptive schedule within the windows exists if and only if no machine is overloaded on any time interval \citep{Horn1974}, and such a schedule is then obtained by Jackson's preemptive rule. Their CP solver, built on this observation, outperforms all previous methods on the pJSSP.

To the best of our knowledge, the fragment-free view has not been exploited for the flexible problem, where it raises two difficulties. First, the set of operations processed by a machine is not known in advance: the overload condition must be enforced on \emph{optional} operations, and its propagation should also remove machines from the alternatives of an operation.
Second, the strength of the approach on the pJSSP relies on a dedicated solver; for the flexible problem, the assignment decisions call for the global reasoning, search strategies and large-neighborhood machinery of a general-purpose CP engine. This work resolves both difficulties and the main contributions are as follows.

\begin{enumerate}
  \item We propose a fragment-free CP model of the pFJSSP, which extends the window model of \citet{Juvin2023CP} to flexible machines. The model has no fragment variable and no time discretization; every solution is turned into a fragment schedule in $O(n \log n)$ time per machine.
  \item We extend the \texttt{PreemptiveNoOverlap} constraint of \citet{Juvin2023CP}, defined for operations whose machine is fixed, to optional intervals, and implement it as a global constraint of IBM ILOG CP Optimizer \citep{Laborie2018}. Optional operations take part in preemptive edge-finding as candidates only, as in the non-preemptive extension of \citet{Vilim2005}, and an absence rule removes infeasible machine alternatives.
  The algorithm runs in $O(n \log n)$ time per call with a $\Theta$-$\Lambda$ tree \citep{Vilim2004} and skips the candidates that cannot tighten a bound. To the best of our knowledge, preemptive edge-finding had not been extended to optional activities.
  \item We combine the CP model with an assignment relaxation based on Carlier's head--tail bounds \citep{Carlier1982}, solved as a mixed-integer program (MILP) with lazily added constraints. The relaxation provides the initial lower bound and, through its solution pool, the machine assignments of the starting solutions. These are completed by the most-work-remaining heuristic of \citet{Juvin2023COR}, implemented on windows.
  \item We report extensive computational experiments on the 276 pFJSSP instances of \citet{Juvin2023COR} and on 78 pJSSP instances. The method improves the best published upper bounds by 2.8\% on average, proves the optimality of 180 instances, and closes 74 of them for the first time. On the pJSSP, we compare it with the solver of \citet{Juvin2023CP} on the same computer and with the same time limit.
\end{enumerate}

The remainder of the paper is organized as follows. Section~\ref{sec:lit} reviews the literature. Section~\ref{sec:problem} defines the problem and presents the theoretical foundations of the fragment-free representation. Section~\ref{sec:model} introduces the CP formulation, and Section~\ref{sec:pno} presents the global constraint. Section~\ref{sec:bounds} describes the lower bounds and initial solutions, while Section~\ref{sec:algorithm} details the complete algorithm. Computational results are reported in Section~\ref{sec:results}. Finally, Section~\ref{sec:conclusion} concludes the paper.

\section{Literature review}
\label{sec:lit}

This section reviews exact approaches for the FJSSP before addressing preemptive shop scheduling, the pFJSSP, and the propagation of resource constraints within CP. We refer to \citet{DauzerePeres2024} for a comprehensive survey of the FJSSP, including its many heuristics.

\subsection{Exact methods for the FJSSP.} The FJSSP generalizes the JSSP, which is strongly NP-hard \citep{Garey1976}. Its classical benchmark sets were introduced by \citet{Brandimarte1993}, \citet{Hurink1994}, \citet{DauzerePeres1997}, \citet{Barnes1996}, \citet{Kacem2002} and \citet{Fattahi2007}. Early exact methods were MILP formulations, which remain limited to small instances; for the distributed FJSSP, \citet{Meng2020} showed that a CP model with interval variables outperforms the MILP formulations. Logic-based Benders decomposition \citep{Hooker2003, Hooker2007} is currently the most effective exact approach: \citet{Naderi2022} assign operations to machines in a MILP master problem and sequence them in JSSP subproblems, and accelerate convergence with cuts derived from the critical paths of the subproblem solutions. Large-neighborhood search and constraint programming have produced most of the best known solutions in recent years \citep[e.g.,][]{Kasapidis2025}.

\subsection{Preemptive shop scheduling.} With preemption, an optimal schedule exists whose preemptions occur at integer times when the data are integral \citep{Baptiste2011}, which justifies discrete-time models. \citet{LePape1998} were the first to solve the pJSSP with CP, using timetabling, disjunctive and edge-finding propagation adapted to preemptive activities \citep[see also][]{Baptiste2001}. \citet{Ebadi2012} proposed MILP models whose size does not depend on the processing times, and \citet{Ebadi2013} a branch-and-bound method that solves instances with up to 50 jobs and 10 machines. \citet{Juvin2023CP} established that, when each operation requires a single resource, the pJSSP can be modelled without representing fragments: with start and end times only, the preemptive no-overlap constraint is equivalent to an overload check, preemptive edge-finding achieves bounds consistency \citep{Baptiste2001}, and the constraint is monotonic, so that a job's windows may be taken back to back. Their solver, implemented in Mistral \citep{Hebrard2008}, proves the optimality of 51 of 63 benchmark instances within one hour. \citet{Terrien2026} extend this approach to maximum workload constraints.

\subsection{The preemptive flexible job shop.} Studies of the pFJSSP are scarce. \citet{Zhang2016} propose metaheuristics for an industrial variant with flexible workdays and overlapping operations. The only exact methods are those of \citet{Juvin2023COR}, described above; their LBBD (denoted LBBD\_B\&B below) relaxes the subproblem in the master problem with machine loads and head--tail constraints, and solves the pJSSP subproblem by branch and bound. Their benchmark of 276 instances, built from the classical FJSSP sets, serves as the problem reference and is used in this paper.

\subsection{Propagation of resource constraints.} Edge-finding was introduced for the one-machine problem by \citet{Carlier1989, Carlier1994} and formalized in CP by \citet{Nuijten1994}; \citet{Baptiste2001} give its preemptive version, which bounds end times instead of start times. \citet{Vilim2004} introduced the $\Theta$-$\Lambda$ tree, which performs overload checking and edge-finding on a unary resource in $O(n \log n)$ time, and \citet{Vilim2005} extended these algorithms to optional activities, whose presence on the resource is not yet decided; such activities are central to flexible scheduling in CP Optimizer, where the machine choice is expressed by an alternative over optional intervals \citep{Laborie2018}. To the best of our knowledge, no preemptive counterpart for optional activities exists; Section~\ref{sec:pno} fills this gap.

\section{Problem definition and fragment-free schedules}
\label{sec:problem}

\subsection{Problem statement}
Let $J$ denote the set of jobs and $M$ the set of machines \cite{Brucker1990}. Each job $j \in J$ consists of an ordered sequence of operations $o_{j,1}, \dots, o_{j,n_j}$ that must be processed in this sequence. The set $O$ represents all operations, while $\text{pred}(o)$ and $\text{succ}(o)$ denote the immediate job predecessor and successor of operation $o$, respectively. 

Each operation $o \in O$ must be assigned to one machine from its set of eligible machines $M_o \subseteq M$, with an integer processing time $p_{o,m}$ on machine $m \in M_o$. We denote by $E_m = \{o \in O : m \in M_o\}$ the set of operations eligible for processing on machine $m$, and by $p^{\min}_o = \min_{m \in M_o} p_{o,m}$ the minimum processing requirement of $o$ across all candidate machines. Each machine can process at most one operation at any given time.

Operations are preemptive: an operation may be interrupted and resumed later, but strictly on the machine to which it was assigned (no inter-machine migration is permitted). A schedule specifies a machine assignment $\mu(o) \in M_o$ for each operation $o \in O$ and allocates a set of disjoint time intervals (\emph{fragments}) whose total duration equals $p_{o,\mu(o)}$. The objective is to minimize the makespan $C_{\max}$, defined as the completion time of the final operation.

The \emph{head} $h(o)$ and \emph{tail} $q(o)$ of an operation $o$ represent the minimum precedence-based processing requirements before and after $o$ within its job:
\begin{equation}
h(o) = \sum_{o' \text{ preceding } o} p^{\min}_{o'}, \qquad 
q(o) = \sum_{o' \text{ succeeding } o} p^{\min}_{o'}.
\end{equation}
Let $S_o$ and $C_o$ denote the start time (start of first fragment) and completion time (end of last fragment) of operation $o$ in a schedule. In every feasible schedule, job precedence dictates that:
\begin{equation}
S_o \ge h(o) \quad \text{and} \quad C_o \le C_{\max} - q(o) \quad \forall o \in O.
\end{equation}

\subsection{Windows, Horn's condition and Jackson's schedule}
\label{sec:windows}
The fragment-free representation rests on three classical scheduling results. For a machine assignment $\mu(o) \in M_o$, the \emph{window} of operation $o$, denoted $[s_o, e_o)$, is defined as any time interval containing all processing fragments of $o$ that lies between the completion of $\text{pred}(o)$ and the start of $\text{succ}(o)$. 

The following proposition establishes necessary and sufficient conditions for the existence of a feasible preemptive schedule on a given machine assignment.

\begin{proposition}[{\citealp{Horn1974, Carlier1982}; \citealp[Theorem~6]{Juvin2023CP}}] \label{prop:horn} Let the machine assignment $\mu$ and operation windows $[s_o, e_o)$ be fixed. The operations assigned to machine $m$ admit a feasible preemptive schedule within their windows if and only if, for every pair of operations $(a, b)$ assigned to machine $m$ with $s_a \le s_b$,
\begin{equation}
\sum_{i : \mu(i)=m,\, s_i \ge s_a,\, e_i \le e_b} p_{i,m} \le e_b - s_a. \label{eq:horn_pairs}
\end{equation} \end{proposition}

\begin{proof}
Condition~\eqref{eq:horn_pairs} is necessary because every operation $i$ whose window is fully contained in the interval $[s_a, e_b)$ must execute entirely within $[s_a, e_b)$. Thus, the total workload of such operations cannot exceed the available machine duration $e_b - s_a$. 

While Horn's original condition considers all non-empty subsets of operations, the pairwise form in~\eqref{eq:horn_pairs} with $s_a \le s_b$ is equivalent. Indeed, consider any arbitrary interval $[s_a, e_{b'})$ containing at least one operation window, and let $b$ be an operation among them with the latest window end ($e_b \le e_{b'}$). Because the window of $b$ is fully contained within $[s_a, e_{b'})$, its start time satisfies $s_b \ge s_a$. Consequently, the pair $(a, b)$ incorporates the exact same set of fully contained operations on the left-hand side of~\eqref{eq:horn_pairs}, while providing an available capacity $e_b - s_a \le e_{b'} - s_a$ that is at least as restrictive. Hence, restricting attention to pairs with $s_a \le s_b$ entails no loss of generality. (Pairs with $e_b < s_a$ yield the vacuous inequality $0 \le e_b - s_a$.)

Sufficiency is constructive: given valid windows, Jackson's preemptive earliest-deadline-first rule \cite{Jackson1955} generates a feasible schedule by executing, at any time, the released operation with the earliest window end deadline $e_o$, preempting it whenever a more urgent operation is released. If~\eqref{eq:horn_pairs} holds, Jackson's schedule satisfies all deadlines and constructs the fragments in $O(n \log n)$ time per machine.
\end{proof}

\begin{proposition}[\citealp{Juvin2023CP}, Lemma~14] \label{prop:monotone} If the operations assigned to a machine admit a feasible preemptive schedule within their windows, they also admit a feasible preemptive schedule within any set of enlarged windows.
\end{proposition}

By Proposition~\ref{prop:monotone}, enlarging windows preserves machine feasibility. Therefore, starting from any valid schedule, windows of consecutive operations in a job can be expanded until they become contiguous, partitioning the job horizon $[0, C_{\max})$ without explicitly tracking processing fragments. Conversely, given a machine assignment and contiguous windows satisfying~\eqref{eq:horn_pairs}, Jackson's preemptive rule recovers an optimal fragment-level schedule.

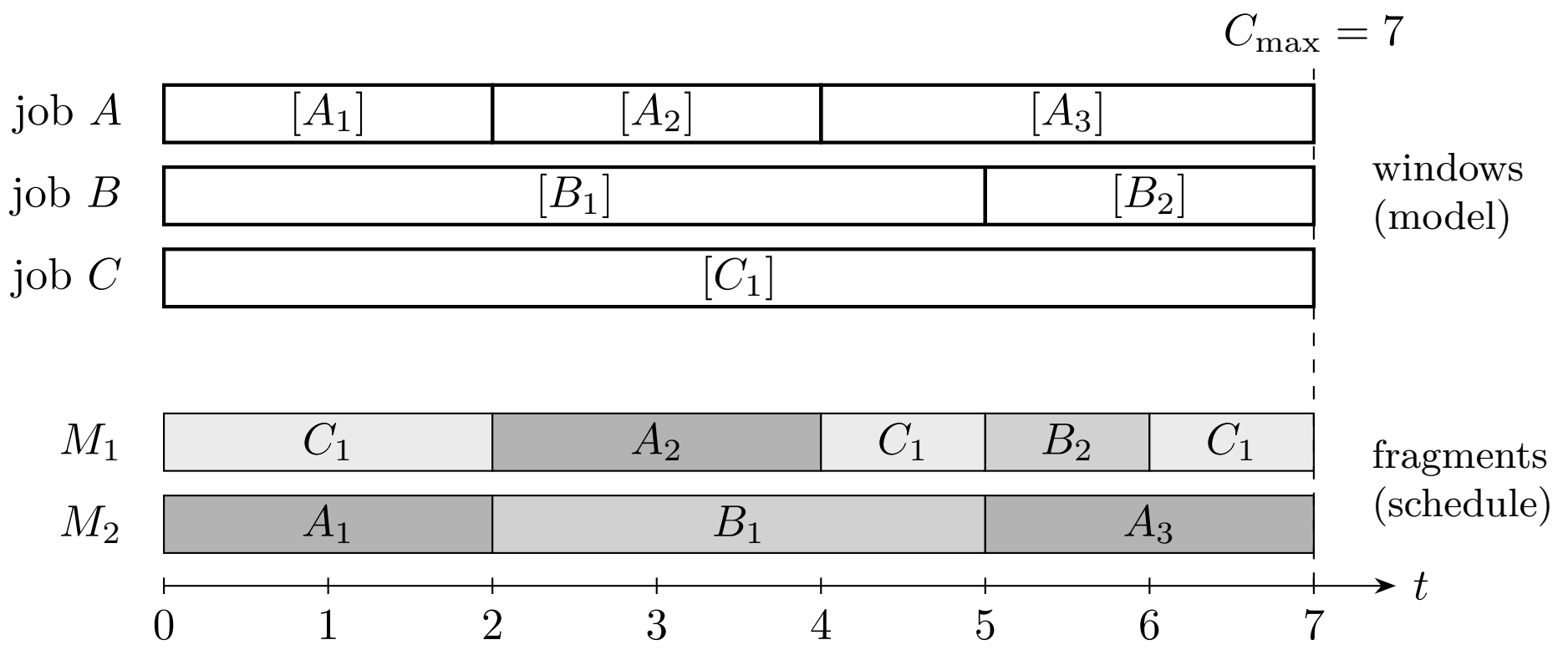
\begin{figure}[htbp]
\centering
\begin{tikzpicture}[x=13mm, y=6.5mm, font=\small]
  \def\win#1#2#3#4{\draw[thick] (#1,#2) rectangle (#3,#2+0.7); \node at ({(#1+#3)/2}, #2+0.35) {#4};}
  \def\frag#1#2#3#4#5{\filldraw[draw=black, fill=#5] (#1,#2) rectangle (#3,#2+0.7); \node at ({(#1+#3)/2}, #2+0.35) {#4};}
  \node[anchor=east] at (-0.15, 6.35) {job $A$};
  \node[anchor=east] at (-0.15, 5.35) {job $B$};
  \node[anchor=east] at (-0.15, 4.35) {job $C$};
  \win{0}{6}{2}{$[A_1]$} \win{2}{6}{4}{$[A_2]$} \win{4}{6}{7}{$[A_3]$}
  \win{0}{5}{5}{$[B_1]$} \win{5}{5}{7}{$[B_2]$}
  \win{0}{4}{7}{$[C_1]$}
  \node[anchor=west, align=left, font=\footnotesize] at (7.25, 5.35) {windows\\(model)};
  \node[anchor=east] at (-0.15, 2.35) {$M_1$};
  \node[anchor=east] at (-0.15, 1.35) {$M_2$};
  \frag{0}{2}{2}{$C_1$}{black!8} \frag{2}{2}{4}{$A_2$}{black!30} \frag{4}{2}{5}{$C_1$}{black!8}
  \frag{5}{2}{6}{$B_2$}{black!18} \frag{6}{2}{7}{$C_1$}{black!8}
  \frag{0}{1}{2}{$A_1$}{black!30} \frag{2}{1}{5}{$B_1$}{black!18} \frag{5}{1}{7}{$A_3$}{black!30}
  \node[anchor=west, align=left, font=\footnotesize] at (7.25, 1.85) {fragments\\(schedule)};
  \draw[-{Stealth}] (0,0.6) -- (7.5,0.6) node[right] {$t$};
  \foreach \t in {0,...,7} {\draw (\t,0.5) -- (\t,0.7); \node[below] at (\t,0.5) {\t};}
  \draw[dashed] (7,0.6) -- (7,6.9) node[above] {$\cmax = 7$};
\end{tikzpicture}
\caption{Windows and fragments. Three jobs on two machines: 
  $A = A_1(M_2, 2) \to A_2(M_1, 2) \to A_3(M_2, 2)$, 
  $B = B_1(M_2, 3) \to B_2(M_1, 1)$, 
  $C = C_1(M_1, 4)$, with the machine assignment and processing time of each operation. 
  \textbf{Top:} the windows of the model, which partition the horizon $[0, C_{\max})$ of each job. 
  \textbf{Bottom:} a preemptive schedule whose fragments lie in these windows; $C_1$ is processed in three fragments. 
  Horn's condition holds on both machines (e.g., on $M_2$, $A_1$ and $B_1$ require 5 units in $[0, 5)$).}
\label{fig:windows}
\end{figure}

Figure~\ref{fig:windows} illustrates this window-based representation.

\section{A fragment-free constraint programming model}
\label{sec:model}
The model uses the interval variables of CP Optimizer \citep{Laborie2018}. For each operation $o$, the interval variable $W_o$ represents its window. Its length is at least the minimum processing requirement, $\mathrm{size}(W_o)\ge p^{\min}_o$. For every eligible machine $m\in M_o$, the optional interval variable $A_{o,m}$ represents the assignment of operation $o$ to machine $m$; it is present if and only if machine $m$ is selected, and then satisfies $\mathrm{size}(A_{o,m})\ge p_{o,m}$. The makespan is denoted by the integer variable $\cmax$.

The CP model is as follows:
\begin{align}
  \min\ & \cmax \label{eq:obj}\\
  \text{s.t.}\ & \texttt{alternative}\big(W_o, \{A_{o,m}\}_{m \in M_o}\big) && \forall o \in O \label{eq:alt}\\
  & \texttt{endAtStart}\big(W_o, W_{\mathrm{succ}(o)}\big) && \forall o \in O,\ o \ne o_{j,n_j} \label{eq:chain}\\
  & \texttt{startOf}(W_{o_{j,1}}) = 0, \quad \texttt{endOf}(W_{o_{j,n_j}}) = \cmax && \forall j \in J \label{eq:ends}\\
  & \textstyle\sum_{o \in E_m} p_{o,m}\, \texttt{presenceOf}(A_{o,m}) \le \cmax && \forall m \in M \label{eq:load}\\
  & \texttt{PreemptiveNoOverlap}\big(\{A_{o,m}\}_{o \in E_m},\ \{p_{o,m}\}_{o \in E_m}\big) && \forall m \in M.
    \label{eq:pno}
\end{align}

Constraints~\eqref{eq:alt} simultaneously perform machine assignment and link the selected assignment interval to the corresponding operation window. Constraints~\eqref{eq:chain}--\eqref{eq:ends} enforce the technological order within each job by partitioning the interval $[0,\cmax)$ into consecutive operation windows. By Proposition~\ref{prop:monotone}, restricting attention to contiguous windows is without loss of generality. Minimizing $\cmax$ therefore compresses all job paths simultaneously. Constraints~\eqref{eq:load} are redundant machine-capacity inequalities. They correspond to the workload condition over the entire planning horizon: the total processing assigned to a machine cannot exceed $\cmax$. Although implied by \eqref{eq:pno}, they provide useful propagation. Constraints~\eqref{eq:pno} enforce the machine-feasibility conditions of Proposition~\ref{prop:horn}. For every machine, they ensure that the workload contained in any interval induced by operation windows does not exceed the available processing capacity. Their implementation is discussed in Section~\ref{sec:pno}. Finally, the domain of $\cmax$ is $[\mathrm{LB},\mathrm{UB}]$, where LB is the lower bound introduced in Section~\ref{sec:milp} and UB is the makespan of the initial solution.

With constraint~\eqref{eq:pno}, the model is exact. Any feasible solution satisfies the conditions of Proposition~\ref{prop:horn} on every machine and therefore admits a feasible preemptive schedule with the same makespan, which can be constructed by Jackson's rule. Consequently, proving infeasibility or optimality of the CP model proves it for the original pFJSSP. Without~\eqref{eq:pno}, the formulation becomes a relaxation in which the windows satisfy job-precedence constraints but may violate machine-capacity requirements on some time interval. Writing out all inequalities in \eqref{eq:horn_pairs} explicitly is impractical. For a machine with $|E_m|$ eligible operations, there are $O(|E_m|^2)$ candidate intervals, each involving a workload sum over up to $|E_m|$ operations. This results in approximately $1.2\times10^5$ conditional terms for a median benchmark instance and $7.4\times10^7$ for the largest one. An alternative is to generate only violated inequalities during the search and add them as cuts. We use this cut-generation scheme as a baseline in Section~\ref{sec:res-prop}.

\section{The \texttt{PreemptiveNoOverlap} constraint on optional intervals}
\label{sec:pno}

The role of constraint~\eqref{eq:pno} is to enforce machine feasibility. For a fixed machine assignment, feasibility is characterized by Proposition~\ref{prop:horn}: the operations assigned to a machine must admit a preemptive schedule within their windows. \citet{Juvin2023CP} introduced the constraint \texttt{PreemptiveNoOverlap} to represent this requirement and showed that it is equivalent to Horn's workload condition. They further showed that the preemptive edge-finding algorithm of \citet{Baptiste2001} establishes bounds consistency and that, by Proposition~\ref{prop:monotone}, propagation can only tighten windows by decreasing latest start times and increasing earliest end times. In the pFJSSP, however, the set of operations assigned to a machine is not known a priori. For each machine $m$, we therefore post a single \texttt{PreemptiveNoOverlap} constraint over the optional assignment intervals $\{A_{o,m}:o\in E_m\}$. The constraint requires that all operations eventually assigned to machine $m$ admit a feasible preemptive schedule within their current windows. During the search, each operation $k\in E_m$ is characterized by \[ \est_k=\texttt{startMin}(A_{k,m}),\qquad \lct_k=\texttt{endMax}(A_{k,m}),\qquad p_k=p_{k,m}, \] together with its assignment status: either the operation is already assigned to machine $m$ (sure) or such an assignment remains possible (optional). The interval $[\est_k,\lct_k)$ therefore represents the largest window still available to operation $k$ on machine $m$.

\subsection{Propagation rules}

The propagation relies on two key observations. First, machine-feasibility reasoning is monotone with respect to window enlargement. By Proposition~\ref{prop:monotone}, enlarging a window cannot destroy a feasible preemptive schedule. Consequently, propagation can only derive that an operation must finish sufficiently late (increase of an earliest completion bound) or start sufficiently early (decrease of a latest start bound). This differs fundamentally from non-preemptive edge-finding, whose main effect is to postpone operations. Second, the earliest completion time of a set of operations can be computed exactly. For a set $\Omega$ of operations assigned to machine $m$,

\begin{equation}
  \ect(\Omega) = \max_{\emptyset \ne \Omega' \subseteq \Omega} \Big( \min_{k \in \Omega'} \est_k + \sum_{k \in \Omega'} p_k \Big)
  \label{eq:ect}
\end{equation}
is the completion time of the last operation of $\Omega$ in the schedule that processes a released operation whenever
the machine is idle, and no preemptive schedule completes $\Omega$ earlier. With
$\lct(\Omega) = \max_{k \in \Omega} \lct_k$, the necessary condition~\eqref{eq:horn_pairs} on the operations of $\Omega$
reads $\ect(\Omega) \le \lct(\Omega)$. Formula~\eqref{eq:ect} is the one computed by the $\Theta$-$\Lambda$ tree in
the non-preemptive case, where it is only a lower bound.

In the following, $\Theta$ denotes a set of operations that are already assigned to machine $m$, and $i\notin\Theta$ denotes an operation whose assignment may or may not be fixed. Propositions~\ref{prop:overload} and \ref{prop:endlast} are the overload and preemptive edge-finding rules of \citet{Baptiste2001} and \citet{Juvin2023CP}. Their extension to the pFJSSP introduces an additional possibility: if applying a rule to a candidate machine would make an operation infeasible, that machine assignment can be eliminated. This yields the absence rule of Proposition~\ref{prop:startfirst}.

\begin{proposition}[{Overload; \citealp{Baptiste2001}; \citealp[Theorem~6]{Juvin2023CP}}]
\label{prop:overload}
If $\ect(\Theta) > \lct(\Theta)$, the constraint is violated. When all intervals are fixed, this test is equivalent
to Horn's condition~\eqref{eq:horn_pairs}.
\end{proposition}

\begin{proposition}[{End-last; \citealp{Baptiste2001}; \citealp[Definition~8]{Juvin2023CP}}]
\label{prop:endlast}
If $\ect(\Theta \cup \{i\}) > \lct(\Theta)$ and $i$ is processed on $m$, then $i$ completes last among
$\Theta \cup \{i\}$, and $\texttt{endOf}(A_{i,m}) \ge \ect(\Theta \cup \{i\})$.
\end{proposition}

\begin{proposition}[Start-first; absence of an optional operation]
\label{prop:startfirst}
The symmetric rule, applied to reversed time ($\est'_k = -\lct_k$, $\lct'_k = -\est_k$), bounds
$\texttt{startOf}(A_{i,m})$ from above. If $i$ is optional and the new end min exceeds its end max, or the new start
max is smaller than its start min, then $i$ cannot be processed on $m$, and $A_{i,m}$ is made absent.
\end{proposition}
\begin{proof}
Reversing time maps preemptive schedules to preemptive schedules and end times to start times, so the start-first
rule is the end-last rule on reversed time. The absence rule is new, since the rules of
\citet{Baptiste2001} and \citet{Juvin2023CP} concern operations whose machine is fixed: both rules assume that $i$
is processed on $m$, and if this assumption empties the domain of $A_{i,m}$, then $i$ is not processed on $m$.
\end{proof}

The absence rule is the main mechanism through which machine flexibility and machine-feasibility reasoning interact. When a machine assignment is ruled out, constraint~\eqref{eq:alt} restricts the operation to its remaining eligible machines; if none remain, infeasibility is detected. Importantly, propagation never increases an earliest start time. Even if an operation cannot be processed before a certain date on a given machine, its window may start earlier because only the processing fragments, not the entire window, must lie on that machine. Raising the earliest start time would therefore be a dominance argument rather than a logical consequence of machine feasibility. Proposition~\ref{prop:monotone} further shows that such a reduction cannot strengthen the machine-feasibility tests.

\emph{Example.} Consider one machine and four operations: $a$ and $b$, sure, with $[\est, \lct) = [0, 4)$ and
$p = 2$; $i$, sure, with $[1, 10)$ and $p = 3$; $j$, optional, with $[0, 6)$ and $p = 3$. For $\Theta = \{a, b\}$,
$\ect(\Theta) = 4 = \lct(\Theta)$: there is no overload. For $i$, $\ect(\{a, b, i\}) = \max(0 + 7,\ 1 + 3) = 7 > 4$,
so $i$ completes last and its end min increases from 4 to 7. For $j$, $\ect(\{a, b, j\}) = 7 > 4$ gives an end min of
7, larger than its end max of 6: $j$ cannot be processed on this machine.

\subsection{Implementation with a \texorpdfstring{$\Theta$-$\Lambda$}{Theta-Lambda} tree}

Applying Propositions~\ref{prop:overload}--\ref{prop:startfirst} to all sets of operations defined by task
intervals takes $O(n^3)$ time per call, where $n$ is the number of non-absent operations of the machine. We use the
$\Theta$-$\Lambda$ tree of \citet{Vilim2004}, with optional activities handled as in \citet{Vilim2005}: the data
structure and its recurrences are theirs. The tree computes~\eqref{eq:ect}, so the detection test is the standard one; only the update differs, since the end min of the candidate $i$ becomes $\ect(\Theta \cup \{i\})$, instead of its start min becoming $\ect(\Theta)$. This preemptive version of the tree, for operations whose machine is fixed,
is the one implemented by \citet{Juvin2023CP} in Mistral. Our additions are the optional candidates and the absence rule, the test of Lemma~\ref{lem:useful}, which discards candidates that cannot tighten a bound, and the implementation choices described in Section~\ref{sec:cpo-integration}. The leaves are the operations sorted by est. Each operation is \emph{white} (in $\Theta$), \emph{gray} (a candidate $i$, in
$\Lambda$) or removed, and each node $v$ stores, over the leaves of its subtree,
\[
  \Sigma P_v = \textstyle\sum_{\text{white}} p, \qquad \ect_v = \ect(\text{white leaves}),
\]
and $\overline{\Sigma P}_v$ and $\overline{\ect}_v$, the maxima of the same quantities when at most one gray leaf is counted as white.
They are computed from the children $l$ (earlier est) and $r$ of $v$:
\begin{align*}
  \Sigma P_v &= \Sigma P_l + \Sigma P_r, &
  \ect_v &= \max\big(\ect_r,\ \ect_l + \Sigma P_r\big),\\
  \overline{\Sigma P}_v &= \max\big(\overline{\Sigma P}_l + \Sigma P_r,\ \Sigma P_l + \overline{\Sigma P}_r\big), &
  \overline{\ect}_v &= \max\big(\overline{\ect}_r,\ \ect_l + \overline{\Sigma P}_r,\ \overline{\ect}_l + \Sigma P_r\big).
\end{align*}
A white leaf has $\Sigma P = \overline{\Sigma P} = p$ and $\ect = \overline{\ect} = \est + p$; a gray leaf has
$\Sigma P = 0$, $\ect = -\infty$, $\overline{\Sigma P} = p$ and $\overline{\ect} = \est + p$. The gray operation
responsible for $\overline{\ect}_{\text{root}}$ is found by a descent from the root that follows the term achieving the maximum at each node. Algorithm~\ref{alg:endlast} gives the end-last pass; the start-first pass is identical on reversed time.

\begin{algorithm}[htbp]
\caption{End-last pass of \texttt{PreemptiveNoOverlap} on one machine}
\label{alg:endlast}
\begin{algorithmic}[1]
\State $S \gets$ sure operations; build the tree with $S$ white and every optional operation $i$ gray if
  $\max(\ect(S), \est_i) + p_i > \mathrm{endMin}(i)$ \Comment{Lemma~\ref{lem:useful}}
\For{each $j \in S$, by non-increasing $\lct_j$} \Comment{white operations $\Theta$: $j$ and the ones not yet visited}
  \If{$\ect_{\text{root}} > \lct_j$} \textbf{fail} \Comment{Proposition~\ref{prop:overload}}
  \EndIf
  \While{$\overline{\ect}_{\text{root}} > \lct_j$}
    \State $i \gets$ gray operation responsible for $\overline{\ect}_{\text{root}}$
    \State $\mathrm{endMin}(i) \gets \max\big(\mathrm{endMin}(i), \overline{\ect}_{\text{root}}\big)$
      \Comment{Proposition~\ref{prop:endlast}; absence if above $\lct_i$}
    \State remove $i$ from the tree
  \EndWhile
  \If{$\ect_{\text{root}} > \mathrm{endMin}(j)$} make $j$ gray \textbf{else} remove $j$
    \Comment{Lemma~\ref{lem:useful}}
  \EndIf
\EndFor
\end{algorithmic}
\end{algorithm}

Many candidate operations cannot yield a stronger bound than the one already
known. Evaluating such candidates is therefore unnecessary. The following
lemma provides a simple dominance test that identifies these cases before they
enter the $\Theta$-$\Lambda$ tree. In instrumented runs, a propagation pass
evaluates between 12 and 32 candidates on average, fewer than one of which
ultimately improves an end-min bound.

\begin{lemma}
\label{lem:useful}
Let $S$ be the set of sure operations, $\Theta \subseteq S$, and
$i \notin \Theta$. Then
\[
  \ect(\Theta \cup \{i\})
  \le
  \max(\ect(S),\est_i)+p_i .
\]
Moreover, if $j$ is made gray when the set of white operations is
$W \ni j$, every later detection of $j$ has value at most $\ect(W)$.
\end{lemma}

\begin{proof}
Let $\Omega' \subseteq \Theta \cup \{i\}$ attain the maximum in
\eqref{eq:ect}. If $\Omega'=\{i\}$, the value is $\est_i+p_i$.
Otherwise,
$\Omega''=\Omega'\setminus\{i\}$
is a non-empty subset of $S$, and
\[
  \min_{k\in\Omega'}\est_k + \sum_{k\in\Omega'} p_k
  \le
  \min_{k\in\Omega''}\est_k + \sum_{k\in\Omega''} p_k + p_i
  \le
  \ect(S)+p_i .
\]
The first statement follows. The second follows from the monotonicity of
$\ect$: later detections of $j$ involve sets
$\Theta \subseteq W\setminus\{j\}$.
\end{proof}

\begin{theorem}
\label{thm:complexity}
Algorithm~\ref{alg:endlast} applies the end-last rule of
Proposition~\ref{prop:endlast} to the sets
\[
  \Theta_j=
  \{k\in S:\lct_k\le \lct_j\},
  \qquad j\in S,
\]
detects every overload of Proposition~\ref{prop:overload} on these sets, and
runs in $O(n\log n)$ time. Every operation receives the strongest end-min
update obtainable from these sets, and the operations discarded by
Lemma~\ref{lem:useful} are precisely those whose current end-min bound cannot
be improved.
\end{theorem}

\begin{proof}
The loop visits the sets $\Theta_j$ by decreasing size. As shown by
\citet{Vilim2004}, an overload of some set of sure operations implies an
overload of one of the sets $\Theta_j$. Each operation is grayed and removed at
most once, and each such event, as well as the descent that identifies the
responsible operation, costs $O(\log n)$. After sorting, the tree is built
bottom-up in $O(n)$ time.

A gray operation is removed at its first detection, which occurs with the
largest set $\Theta_j$ and therefore yields its strongest end-min update. Since
gray operations do not belong to $\Theta$, discarding some of them affects
neither the overload tests nor the propagation of other candidates. The final
statement follows directly from Lemma~\ref{lem:useful}.
\end{proof}

\subsection{Integration in CP Optimizer}
\label{sec:cpo-integration}
Propagation by \texttt{PreemptiveNoOverlap} interacts closely with the machine assignment, job-precedence and makespan constraints. Whenever the feasible window of an operation changes, the propagator re-evaluates machine feasibility and may derive stronger end-min bounds, stronger start-max bounds, or the elimination of some machine assignments. These reductions are then propagated through the alternative constraints~\eqref{eq:alt}, the job chains \eqref{eq:chain}--\eqref{eq:ends}, and the objective. For example, increasing the end minimum of $A_{i,m}$ enlarges the window $W_i$, which may in turn shift the windows of successor operations and increase the lower bound on $\cmax$. Propagation continues until a fixpoint is reached. At the implementation level, \texttt{PreemptiveNoOverlap} is realized as a dedicated engine-level constraint of CP Optimizer. It is triggered by changes to interval bounds or presence decisions and executes once per propagation round, applying the overload, end-last and start-first rules described above. Three implementation choices reduce the cost of a propagation call. First, the sorting of operations by $\est$ reuses the order from the previous call on the same machine through insertion sort, since only a few bounds typically change between successive calls. Second, the operation responsible for $\overline{\ect}_{\mathrm{root}}$ is identified on demand rather than stored in the tree nodes. Third, the propagator caches its previous inputs and outputs and reuses them whenever the operation windows have not changed, avoiding the reconstruction of the $\Theta$-$\Lambda$ trees. On recorded calls, these optimizations reduce the running time of a propagation call by a factor between 1.6 and 1.9 without affecting the results. On the largest benchmark instances, \texttt{PreemptiveNoOverlap} accounts for approximately 45\% to 50\% of the total search time. 

\section{Lower bound and starting solutions}
\label{sec:bounds}
\subsection{An assignment relaxation with Carlier's bounds} \label{sec:milp} The lower bound is the optimal value of an assignment relaxation, which is the master problem of the LBBD of \citet{Juvin2023COR} without Benders cuts. Let $x_{o,m}\in\{0,1\}$ indicate that operation $o$ is assigned to machine $m$, and let $C\ge0$ denote the makespan: \begin{align} \min\ & C \nonumber\\ \text{s.t.}\ & \textstyle\sum_{m \in M_o} x_{o,m} = 1 && \forall o \in O \label{eq:assign}\\ & \textstyle C \ge \sum_{o \in \text{job } j} \sum_{m \in M_o} p_{o,m}\, x_{o,m} && \forall j \in J \label{eq:job}\\ & \textstyle C \ge h + \sum_{o \in S(m,h,q)} p_{o,m}\, x_{o,m} + q && \forall m \in M,\ (h,q)\in T_m, \label{eq:carlier} \end{align} where \[ T_m=\{(\min(h(a),h(b)),\ \min(q(a),q(b))) : a,b\in E_m\} \] and \[ S(m,h,q)=\{o\in E_m : h(o)\ge h,\ q(o)\ge q\}. \] 

Constraints~\eqref{eq:assign} assign each operation to one eligible machine. Constraints~\eqref{eq:job} are job-based lower bounds. Constraints \eqref{eq:carlier} are Carlier's single-machine head--tail bounds \cite{Carlier1982}: the operations of $S(m,h,q)$ assigned to machine $m$ cannot start before $h$ and must complete by $C-q$, whether preemption is allowed or not. The formulation is thus a valid relaxation of the pFJSSP. 

The model may contain up to $\sum_m |E_m|^2$ constraints~\eqref{eq:carlier}, most of which are never binding. 
Furthermore, the relaxation becomes computationally demanding when the bound arises from balancing flexible workload with similar processing times. We therefore add constraints~\eqref{eq:carlier} lazily in rounds: the model is initially solved with \eqref{eq:assign}, \eqref{eq:job}, and basic machine workload bounds ($h = q = 0$). Violated constraints \eqref{eq:carlier} are identified and added, and the model is re-solved, warm-started from the current solution, until an optimal round violates no constraint. We found this multi-round iterative scheme to yield stronger bounds within the time limit than a single branch-and-bound tree with lazy constraint or user cut callbacks. Because each round represents a valid relaxation, its rounded-up objective value provides a valid lower bound even if interrupted by the time limit.

\subsection{Starting solutions} \label{sec:start} 
Two complementary procedures construct starting primal solutions for the CP search.

\subsubsection{Preemptive EDD on Relaxation Assignments.} Given a machine assignment $x$ and operation due dates $d(o)$, preemptive earliest-due-date (EDD) list scheduling releases each operation upon completion of its job predecessor and assigns machine priority to released operations with the earliest due date (Jackson's rule). We apply EDD to every integer assignment in the MILP solution pool (the incumbent and alternative integer solutions from the final round) setting:
\begin{equation}
d(o) = \text{LB} - q(o),
\end{equation}
which represents the latest completion time of operation $o$ consistent with the lower bound $\text{LB}$. While the MILP relaxation balances capacity without considering intra-machine sequence, EDD recovers feasible sequencing.

\subsubsection{Window Greedy.} The second construction simultaneously determines machine assignments and schedule placement. It follows the assignment logic of Juvin et al. \citet{Juvin2023COR} by prioritizing the job with the most work remaining (MWKR), a classical dispatching rule well suited for makespan minimization \citep{Chang1996}, on the machine that completes it earliest. However, rather than fixing discrete processing fragments, it constructs operation windows directly, preserving internal resequencing flexibility. At each step, the next operation of every job is evaluated on each eligible machine with candidate window $[r, e)$, where $r$ is the window end of its job predecessor and $e$ is the earliest completion time preserving machine feasibility (determined via binary search over Jackson's schedule). The job maximizing total remaining work ($\sum_{o' \text{ unscheduled}} p^{\min}_{o'}$) is committed to the machine attaining the earliest window end $e$. Once all windows are placed, Jackson's rule constructs the final fragment schedule. The greedy construction executes in at most 0.06~s across all benchmark instances. 
The best schedule generated across EDD pool evaluation and the window greedy heuristic serves as the warm-start incumbent for the CP search.

\section{The complete algorithm} \label{sec:algorithm} 
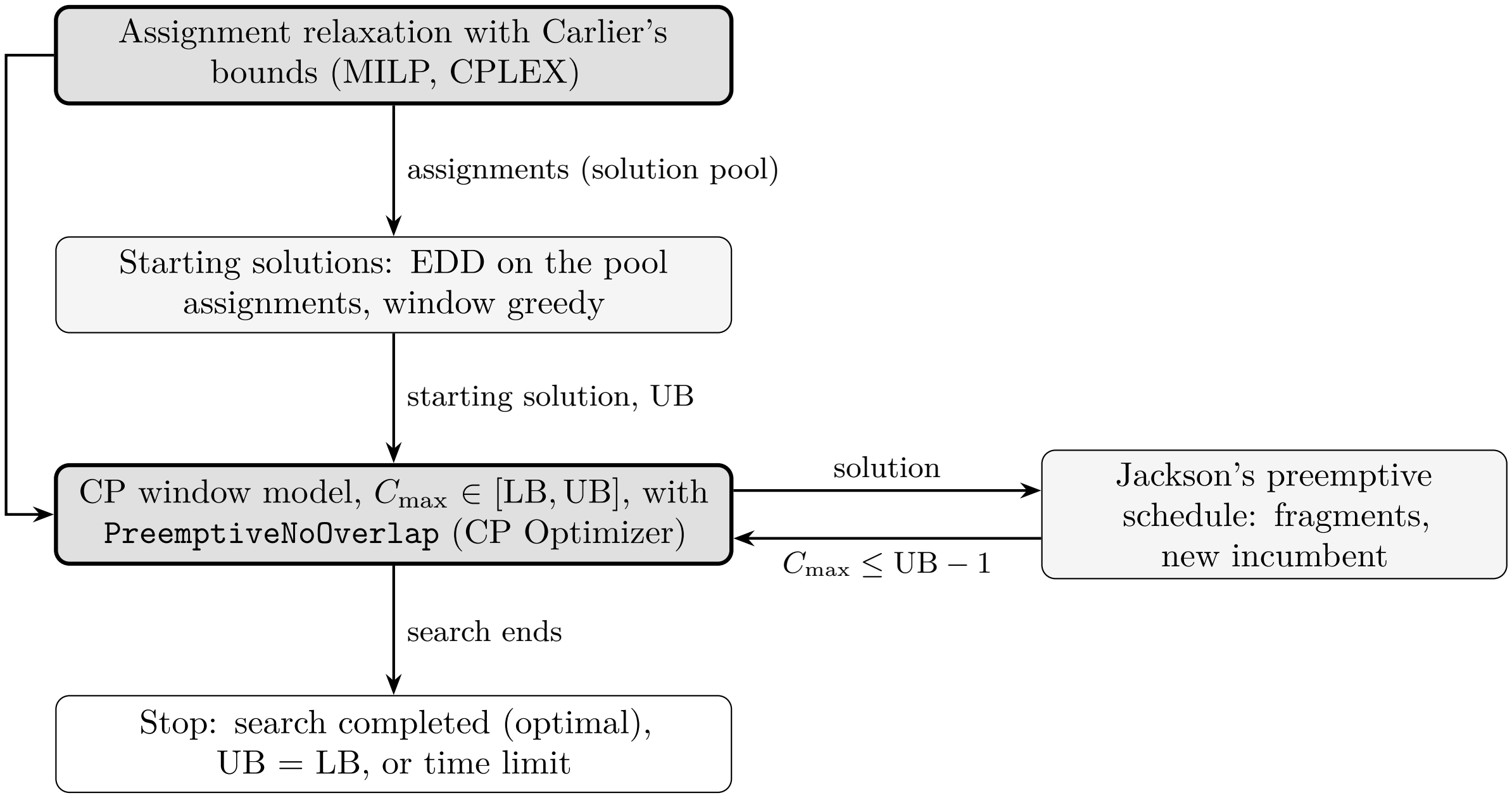
\begin{figure}[htbp]
\centering
\begin{tikzpicture}[font=\small,
  box/.style={draw, rounded corners, align=center, minimum height=10mm, text width=68mm, fill=black!4},
  opt/.style={box, fill=black!12, very thick},
  arr/.style={-{Stealth[length=2.2mm]}, thick},
  lab/.style={font=\footnotesize}]
  \node[opt] (milp) at (0, 0) {Assignment relaxation with Carlier's\\bounds (MILP, CPLEX)};
  \node[box] (start) at (0, -2.4) {Starting solutions: EDD on the pool\\assignments, window greedy};
  \node[opt] (cp) at (0, -4.8) {CP window model, $\cmax \in [\mathrm{LB}, \mathrm{UB}]$, with\\
    \texttt{PreemptiveNoOverlap} (CP Optimizer)};
  \node[box, text width=46mm] (jack) at (9.2, -4.8) {Jackson's preemptive\\schedule: fragments,\\new incumbent};
  \node[box, fill=white] (stop) at (0, -7.2) {Stop: search completed (optimal),\\UB $=$ LB, or time limit};
  \draw[arr] (milp) -- node[lab, right] {assignments (solution pool)} (start);
  \draw[arr] (start) -- node[lab, right] {starting solution, UB} (cp);
  \draw[arr] (milp.west) -- ++(-0.5, 0) |- node[lab, left, pos=0.25] {LB} (cp.west);
  \draw[arr] ([yshift=2.5mm]cp.east) -- node[lab, above] {solution} ([yshift=2.5mm]jack.west);
  \draw[arr] ([yshift=-2.5mm]jack.west) -- node[lab, below] {$\cmax \le \mathrm{UB} - 1$} ([yshift=-2.5mm]cp.east);
  \draw[arr] (cp) -- node[lab, right] {search ends} (stop);
\end{tikzpicture}
\caption{Overview of \cpw. The two optimization models (shaded) are solved by CPLEX and CP Optimizer; fragments are
never variables, and are built by Jackson's preemptive schedule from the machines and windows of each solution.}
\label{fig:pipeline}
\end{figure}

Algorithm~\ref{alg:cpw}, denoted \cpw{}, combines the assignment relaxation and the window-based CP model. Figure~\ref{fig:pipeline} summarizes the workflow. The assignment relaxation is solved first, using 10\% of the time limit. Its objective value gives the lower bound LB, while the best starting solution provides the upper bound UB and the starting point of the CP search. CP Optimizer then searches the model of Section~\ref{sec:model} with its Restart search \citep{Laborie2018}, which combines restarts, large neighborhood search and failure-directed search \citep{Vilim2015}. Every solution specifies machines and windows; the corresponding fragment schedule is reconstructed by Jackson's rule and becomes the incumbent. The algorithm stops when the search is complete, when $\mathrm{UB}=\mathrm{LB}$, or when the time limit is reached.

\begin{algorithm}[htbp]
\caption{\cpw}
\label{alg:cpw}
\begin{algorithmic}[1]
\State $(\mathrm{LB}, X) \gets$ lazy Carlier MILP: lower bound and pool of assignments \Comment{Section~\ref{sec:milp}}
\State $\sigma^* \gets$ best of $\{\mathrm{EDD}(x, d = \mathrm{LB} - q) : x \in X\} \cup \{\text{window greedy}\}$;
  $\mathrm{UB} \gets \cmax(\sigma^*)$ \Comment{Section~\ref{sec:start}}
\State build the model \eqref{eq:obj}--\eqref{eq:pno} with $\cmax \in [\mathrm{LB}, \mathrm{UB}]$; starting point $\sigma^*$
\While{the search finds a solution $\sigma$ and time remains}
  \State $\sigma^* \gets$ Jackson's schedule of $\sigma$; $\mathrm{UB} \gets \cmax(\sigma)$
    \Comment{feasible by Proposition~\ref{prop:horn}}
\EndWhile
\State \textbf{if} the search is complete \textbf{then} $\mathrm{LB} \gets \mathrm{UB}$ \Comment{$\sigma^*$ is optimal}
\end{algorithmic}
\end{algorithm}

\section{Computational experiments} \label{sec:results} \subsection{Experimental setting}

\cpw{} is implemented in C++17 with IBM ILOG CP Optimizer and CPLEX~22.2.0. All experiments were run on an Apple M4 processor (macOS, native arm64 build), using a single thread per run, four runs in parallel, and a time limit of 600\,s per instance, including the MILP phase. Unless stated otherwise, \cpw{} is used with the configuration of Section~\ref{sec:algorithm} and the same parameter settings on all instances. 
\subsubsection{Instances.} The pFJSSP benchmark of \citet{Juvin2023COR} comprises 276 instances derived from the FJSSP sets of \citet{Barnes1996} (21 instances), \citet{Brandimarte1993} (15), \citet{DauzerePeres1997} (18), \citet{Fattahi2007} (20), \citet{Hurink1994} (edata, rdata and vdata, 66 each) and \citet{Kacem2002} (4). The instances contain 2 to 30 jobs, 2 to 18 machines, and between 1 and 10 eligible machines per operation on average. For the pJSSP, we use the 78 instances considered by \citet{Terrien2026}, including the 63 instances of \citet{Juvin2023CP} from \citet{Fisher1963}, \citet{Lawrence1984}, \citet{Adams1988}, \citet{Applegate1991} and \citet{Storer1992} (swv16--20), together with swv01--15 \citep{Storer1992}, for which no preemptive results have been published. 

\subsubsection{Reference methods.} For the pFJSSP, we compare with the four methods of \citet{Juvin2023COR}: CP Optimizer on their model with and without a warm start, and the two LBBD variants, including LBBD\_B\&B, the best-performing method in their study. Their experiments were run with a one-hour time limit on an Intel Xeon E5-2695 v4 processor at 2.1\,GHz using CPLEX and CP Optimizer. Because the hardware differs, we focus primarily on lower bounds, proofs of optimality and the number of closed instances. For the pJSSP, we reran the solver of \citet{Juvin2023CP} (Mistral, default settings) on the same machine and with the same time limits as \cpw. This comparison is therefore free of hardware effects. Our rerun reproduced all optimal values reported by \citet{Juvin2023CP}. 

\subsubsection{Measures.} For an instance, the gap of an upper bound UB to a reference value $R$ is \[ 100(\mathrm{UB}-R)/R. \] The optimality gap of a method is \[ 100(\mathrm{UB}-\mathrm{LB})/\mathrm{LB}, \] using its own lower and upper bounds. An instance is \emph{closed} if its optimal value is known. It is \emph{newly closed} if its optimal value was not established by the published results, nor by a non-preemptive optimum matching the best published lower bound.

\subsection{Preemptive job shop} \label{sec:res-pjssp} 
Table~\ref{tab:pjssp} compares \cpw{} with Mistral on the 78 pJSSP instances. \begin{table}[htbp] \centering\small \caption{Preemptive job shop: \cpw{} and Mistral \citep{Juvin2023CP} on the same computer, with the same time limits. Better / equal / worse: comparison of the upper bounds of \cpw{} with those of Mistral.} \label{tab:pjssp} 

\begin{tabular}{ll r rr rrr rr}
\toprule
 & & & \multicolumn{2}{c}{Proven optimal} & \multicolumn{3}{c}{UB of CPW vs Mistral} & \multicolumn{2}{c}{Time (s)$^a$} \\
\cmidrule(lr){4-5}\cmidrule(lr){6-8}\cmidrule(lr){9-10}
Limit & Instances & $n$ & Mistral & CPW & Better & Equal & Worse & Mistral & CPW \\
\midrule
60 s & CP 2023 set & 63 & 45 & 60 & 12 & 51 & 0 & 132 & 16.1 \\
 & swv01--15 & 15 & 1 & 4 & 10 & 1 & 4 & 4 & 2.7 \\
 & All & 78 & 46 & 64 & 22 & 52 & 4 & 136 & 18.9 \\
\midrule
600 s & CP 2023 set & 63 & 50 & 61 & 9 & 54 & 0 & 697 & 24.8 \\
 & swv01--15 & 15 & 1 & 6 & 14 & 1 & 0 & 3 & 2.9 \\
 & All & 78 & 51 & 67 & 23 & 55 & 0 & 700 & 27.7 \\
\bottomrule
\end{tabular}

\smallskip \begin{minipage}{0.92\textwidth}\footnotesize $^a$ Total time to the proof on the instances proven by both methods. \end{minipage} \end{table} \cpw{} clearly dominates Mistral on this benchmark. It proves 64 instances within 60\,s and 67 within 600\,s, compared with 46 and 51 for Mistral. Mistral proves no instance that \cpw{} does not prove, and both methods agree on every optimal value. At 600\,s, \cpw{} obtains a better upper bound on 23 instances and is never worse. On the instances proven by both methods, \cpw{} requires a total of 28\,s, against 700\,s for Mistral. Part of the improvement comes from the lower bound used by \cpw{}, which, on the pJSSP, reduces to the best of Carlier's one-machine bounds and the job length bounds. However, the propagation of \texttt{PreemptiveNoOverlap} is also critical. With the same lower bound and starting solution, but replacing the propagator by lazy cuts, \cpw{} proves only 29 instances within 600\,s, fewer than Mistral (Section~\ref{sec:res-prop}). For reference, \citet{Juvin2023CP} reported 51 proofs of optimality on the 63 instances of their benchmark with a one-hour time limit. Our rerun proves 50 instances in 600\,s, confirming that Mistral is a strong baseline. Among the 63 benchmark instances of \citet{Juvin2023CP}, eleven remained open. Within 600\,s, \cpw{} closes nine of them: abz7 (651), la24 (909), la25 (947), la29 (1126), la36 (1252), la38 (1166), la40 (1198), orb03 (973), and orb06 (985). Only abz8 and abz9 remain open. On the additional instances swv01--15, for which no preemptive results have been published, \cpw{} proves the optimality of swv01 (1380), swv02 (1475), swv03 (1372), swv05 (1411), swv13 (3104), and swv14 (2968), and improves the best upper bound obtained by Mistral on 14 of the 15 instances. For swv13 and swv14, the preemptive and non-preemptive optima coincide.

\subsection{Preemptive flexible job shop} \label{sec:res-pfjssp} Table~\ref{tab:pfjssp} summarizes the results on the 276 pFJSSP instances by family. The complete results are available in the online repository of the paper (link provided upon acceptance). \begin{table}[htbp] \centering\small \caption{Preemptive flexible job shop: \cpw{} (600\,s) against the methods of \citet{Juvin2023COR} (3600\,s). UB: mean gap of the upper bound of \cpw{} to that of LBBD\_B\&B and to the best of the four published methods, and number of instances on which it is at most that of LBBD\_B\&B.} \label{tab:pfjssp} 
\resizebox{\textwidth}{!}{
\begin{tabular}{l r rr r rr rr c r}
\toprule
 & & \multicolumn{3}{c}{UB of CPW} & \multicolumn{2}{c}{Proven optimal} & \multicolumn{2}{c}{Mean gap (\%)$^a$} & LB of CPW$^b$ & Newly \\
\cmidrule(lr){3-5}\cmidrule(lr){6-7}\cmidrule(lr){8-9}
Family & $n$ & vs LBBD (\%) & vs best (\%) & $\le$ LBBD & CPW & LBBD & CPW & LBBD & $>$ / $=$ / $<$ & closed \\
\midrule
Barnes & 21 & $-$5.45 & $-$4.64 & 21 & 21 & 0 & 0.00 & 11.20 & 20 / 1 / 0 & 21 \\
Brandimarte & 15 & $-$0.44 & +0.73 & 11 & 8 & 10 & 4.11 & 4.86 & 2 / 13 / 0 & 0 \\
Dauz\`ere-P\'er\`es & 18 & $-$7.34 & $-$6.89 & 18 & 2 & 0 & 4.16 & 12.62 & 3 / 14 / 1 & 0 \\
Fattahi & 20 & +0.18 & +0.18 & 19 & 19 & 17 & 0.70 & 0.83 & 2 / 18 / 0 & 2 \\
Hurink edata & 66 & $-$2.85 & $-$2.49 & 64 & 60 & 26 & 0.62 & 7.24 & 36 / 30 / 0 & 32 \\
Hurink rdata & 66 & $-$5.87 & $-$4.58 & 52 & 30 & 18 & 2.48 & 10.42 & 23 / 41 / 2 & 19 \\
Hurink vdata & 66 & $-$2.62 & $-$1.59 & 55 & 36 & 27 & 0.62 & 3.45 & 0 / 66 / 0 & 0 \\
Kacem & 4 & +0.00 & +0.00 & 4 & 4 & 4 & 0.00 & 0.00 & 0 / 4 / 0 & 0 \\
\midrule
All & 276 & $-$3.62 & $-$2.82 & 244 & 180 & 102 & 1.44 & 7.05 & 86 / 187 / 3 & 74 \\
\bottomrule
\end{tabular}}
\smallskip 
\begin{minipage}{0.97\textwidth}\footnotesize $^a$ Mean optimality gap, each method with its own bounds. $^b$ Lower bound of \cpw{} higher than, equal to, or lower than that of LBBD\_B\&B. \end{minipage} \end{table} \cpw{} substantially improves the state of the art on this benchmark. Its upper bound is at most that of LBBD\_B\&B on 244 of the 276 instances and improves upon it by 3.62\% on average. Relative to the best upper bound reported by any method of \citet{Juvin2023COR}, the average improvement is 2.82\%. The largest gains occur on the Dauz\`ere-P\'er\`es instances ($-7.34\%$), Hurink rdata ($-5.87\%$), and Barnes ($-5.45\%$). The improvement is reflected in both primal and exact results. \cpw{} proves the optimality of 180 instances, compared with 102 for LBBD\_B\&B, and reduces the mean optimality gap from 7.05\% to 1.44\%. Its lower bound is higher than that of LBBD\_B\&B on 86 instances and lower on only three. Most proofs are obtained quickly: 117 instances are solved within one second and 171 within one minute. Combining its bounds with the published results closes 205 instances, including 74 that were previously open: all 21 Barnes instances, 32 Hurink edata instances, 19 Hurink rdata instances, and two Fattahi instances. In particular, \cpw{} proves the optimality of every Barnes instance, whereas no published method had closed this family completely. Despite these results, \cpw{} does not dominate LBBD\_B\&B on every instance. Its upper bound is worse on 32 instances. On 21 of them, however, the difference is only one to three time units above an optimum proven by LBBD\_B\&B; this situation occurs mainly in the Hurink rdata, Hurink vdata, and Brandimarte families, and in 19 cases the best bound is reached within six seconds and remains unchanged thereafter. Larger differences are observed on a small number of instances, including mfjs10 (1274 versus 1204) and v-car5 (4961 versus 4912). Finally, the upper bound of \cpw{} remains above the best known non-preemptive makespan \citep{Kasapidis2025} on 94 instances. Since every non-preemptive schedule is feasible for the preemptive problem, these instances indicate remaining room for improvement on the primal side of the method.

\subsection{Contribution of the components} \label{sec:res-components} 

Table~\ref{tab:components} separates the contribution of the first step of \cpw{} (the MILP bound and the best starting solution) from that of the CP search. The first step is inexpensive: the MILP requires 0.02\,s for the median instance and 7.4\,s on average, while the construction of the starting solutions is negligible. \begin{table}[htbp] \centering\small \caption{Contribution of the components: bounds after the first step (MILP lower bound and best starting solution), and after the CP search. Results are shown with lazy overload cuts instead of the propagator (Section~\ref{sec:res-prop}) and with the propagator, 600\,s.} \label{tab:components} 
\resizebox{\textwidth}{!}
{\begin{tabular}{ll r r r r c}
\toprule
Benchmark & Bounds & UB vs LBBD (\%) & UB $\le$ LBBD & Proven & Mean gap (\%) & LB vs LBBD$^a$ \\
\midrule
pFJSSP (276) & First step: MILP bound and best start & +5.96 & 91 & 35 & 13.46 & 17 / 241 / 18 \\
 & CP search with lazy cuts & $-$2.03 & 204 & 95 & 4.37 & 33 / 236 / 7 \\
 & CP search with the propagator (\cpw) & $-$3.62 & 244 & 180 & 1.44 & 86 / 187 / 3 \\
\midrule
pJSSP (78) & First step: MILP bound and best start & -- & -- & 12 & 20.67 & -- \\
 & CP search with lazy cuts & -- & -- & 29 & 6.64 & -- \\
 & CP search with the propagator (\cpw) & -- & -- & 67 & 0.93 & -- \\
\bottomrule
\end{tabular}}
\smallskip \begin{minipage}{0.97\textwidth}\footnotesize $^a$ Lower bound higher than, equal to, or lower than that of LBBD\_B\&B. \end{minipage} \end{table} The first step alone is not competitive with the state of the art. On the pFJSSP, its best starting solution is on average 5.96\% above the upper bound of LBBD\_B\&B and improves upon it on only 91 instances. Its lower bound is essentially that of LBBD\_B\&B, being higher on 17 instances and lower on 18, which is expected since the MILP corresponds to the relaxation used in the LBBD master problem. Only 35 instances are closed at this stage, when the starting solution matches the lower bound. Most of the improvement therefore comes from the CP search. It improves the starting solution on 239 instances and reduces the average gap to the upper bound of LBBD\_B\&B from $+5.96\%$ to $-3.62\%$. Among the 180 instances proven optimal by \cpw{}, 145 are closed during the CP phase. On 63 of them, the search reaches the MILP bound; on the remaining 82, it proves optimality above that bound after completing the search. The CP search also strengthens the lower bound. It improves the MILP bound on 87 instances, by 5.0\% on average on those instances, and is responsible for nearly all of the cases in which \cpw{} obtains a stronger lower bound than LBBD\_B\&B. Overall, the mean optimality gap decreases from 13.46\% after the first step to 1.44\% at the end of the search. The pJSSP exhibits the same pattern. The first step closes 12 instances, whereas the CP search closes a further 55, reducing the mean optimality gap from 20.67\% to 0.93\%.

\subsection{Contribution of the propagator} \label{sec:res-prop} 
To assess the contribution of the \texttt{PreemptiveNoOverlap} global constraint, we evaluate a baseline variant in which the propagator is replaced by a lazy overload cut scheme. In this variant, the CP model is solved as a relaxation that omits constraint~\eqref{eq:pno}. 

Whenever a candidate solution violates Horn's condition~\eqref{eq:horn_pairs} for a pair of operations $(a, b)$ on machine $m$, it is first subjected to a repair step using preemptive EDD list scheduling under the current machine assignments and due dates set to window ends ($d(o) = e_o$). If the repaired schedule fails to improve upon the incumbent makespan, the search is paused, the violated condition is posted as a linear constraint over the decision variables, and the search is restarted from the incumbent solution.

Let $\pi_i = \text{presenceOf}(A_{i,m})$ denote the presence indicator of operation $i$ on machine $m$, $S_i = \text{startOf}(W_i)$ its window start time, and $F_i = \text{endOf}(W_i)$ its window end time. For an overloaded pair $(a, b)$ with interval bounds $S_a = \text{startOf}(W_a)$ and $F_b = \text{endOf}(W_b)$, the lazy cut posted to the model is formulated as:

\begin{equation}
\sum_{i \in E_m} p_{i,m} \cdot \left[\pi_i \land (S_i \ge S_a) \land (F_i \le F_b)\right] \le \max(0, F_b - S_a),
\label{eq:lazy_cut_formulation}
\end{equation}

where $[\cdot] \in \{0, 1\}$ represents the Iverson bracket (returning 1 if the condition holds and 0 otherwise).

All other components of the algorithm remain identical. Table~\ref{tab:lazy} compares the propagator against this lazy-cut variant across all 276 pFJSSP instances and 78 pJSSP instances. 

As shown in Table~\ref{tab:lazy}, the propagator roughly doubles the number of optimality proofs (180 vs.\ 95 on the pFJSSP, and 67 vs.\ 29 on the pJSSP within 600~s). Two key structural factors explain this performance difference:
\begin{enumerate}
    \item \textbf{Search Restarts vs.\ Continuous Filtering:} The lazy variant discovers violated intervals one at a time and must restart the CP search after posting each cut (averaging 15 restarts and ~1,000 posted constraints per instance, reaching up to 68,000 branches/s overhead).
    \item \textbf{Global Scope:} Lazy cuts only restrict specific intervals that have already proven infeasible in prior candidate solutions. In contrast, \texttt{PreemptiveNoOverlap} filters invalid domain values across all operation subsets dynamically at every search node in $O(n \log n)$ time, allowing FDS and LNS search heuristics to prune subtrees substantially earlier.
\end{enumerate}

\begin{table}[htbp] \centering\small \caption{Propagator against lazy overload cuts, all other components identical (600\,s unless stated): instances proven optimal, mean gap of the upper bound to that of LBBD\_B\&B, mean optimality gap, and comparison of the upper bounds of the two variants.} \label{tab:lazy} 
\resizebox{\textwidth}{!}
{
 \begin{tabular}{l r rr rr rr c}
\toprule
 & & \multicolumn{2}{c}{Proven optimal} & \multicolumn{2}{c}{UB vs LBBD (\%)} & \multicolumn{2}{c}{Mean gap (\%)} & UB of propagator \\
\cmidrule(lr){3-4}\cmidrule(lr){5-6}\cmidrule(lr){7-8}
Benchmark & $n$ & Prop. & Lazy & Prop. & Lazy & Prop. & Lazy & better / equal / worse \\
\midrule
Barnes & 21 & 21 & 0 & $-$5.45 & $-$1.48 & 0.00 & 8.43 & 19 / 2 / 0 \\
Brandimarte & 15 & 8 & 8 & $-$0.44 & $-$0.41 & 4.11 & 4.19 & 2 / 11 / 2 \\
Dauz\`ere-P\'er\`es & 18 & 2 & 0 & $-$7.34 & $-$5.58 & 4.16 & 6.15 & 12 / 1 / 5 \\
Fattahi & 20 & 19 & 17 & +0.18 & +0.82 & 0.70 & 1.83 & 3 / 17 / 0 \\
Hurink edata & 66 & 60 & 24 & $-$2.85 & $-$0.45 & 0.62 & 6.09 & 36 / 28 / 2 \\
Hurink rdata & 66 & 30 & 10 & $-$5.87 & $-$3.86 & 2.48 & 5.39 & 40 / 14 / 12 \\
Hurink vdata & 66 & 36 & 33 & $-$2.62 & $-$2.48 & 0.62 & 0.76 & 22 / 37 / 7 \\
Kacem & 4 & 4 & 3 & +0.00 & +2.50 & 0.00 & 2.50 & 1 / 3 / 0 \\
\midrule
All & 276 & 180 & 95 & $-$3.62 & $-$2.03 & 1.44 & 4.37 & 135 / 113 / 28 \\
\midrule
pJSSP, 60 s & 78 & 64 & 28 & -- & -- & 1.70 & 8.98 & 47 / 30 / 1 \\
pJSSP, 600 s & 78 & 67 & 29 & -- & -- & 0.93 & 6.64 & 47 / 31 / 0 \\
\bottomrule
\end{tabular}
}
\end{table} 

Both variants start from the same lower bound and the same starting solution on every instance. The observed differences therefore come entirely from the propagation mechanism. The propagator roughly doubles the number of proofs: 180 against 95 on the pFJSSP and 67 against 29 on the pJSSP. The lazy variant proves no instance that the propagator does not prove, and never obtains a better upper bound on an instance solved to optimality by the propagator. The effect is particularly visible on solved instances. On instances proven by both variants, the propagator requires 67\,s of CPU time in total, against 3110\,s for the lazy variant (median ratio 4.6 on the pFJSSP). The lower bounds show the same pattern: the CP search improves the MILP bound on 87 instances with the propagator, but on only 27 with lazy cuts, and the lower bound obtained with the propagator is never worse. The largest gains occur on families where optimality can be established within the time limit. On the Barnes instances, the propagator proves all 21 instances, whereas the lazy variant proves none. On Hurink edata, the numbers are 60 and 24, respectively. On the more difficult Dauz\`ere-P\'er\`es and Hurink rdata families, where most instances remain open, the lazy variant obtains a better upper bound on 28 instances, by 1.44\% on average, whereas the propagator is better on 135 instances, by 3.82\% on average. On Brandimarte and Hurink vdata, where most instances are either solved quickly or limited by the lower bound, the two variants are closer. Two factors explain this difference. First, the lazy variant discovers violated conditions one solution at a time and must restart the search to post them. It performs 15 searches per instance on average, posting about 1000 constraints, and up to 62 searches and 6866 constraints on a single instance. Second, these constraints propagate weakly and only on the intervals already found violated. In contrast, \texttt{PreemptiveNoOverlap} reasons on all relevant sets of operations at every search node. As a result, each branch is considerably cheaper: on instances solved by neither variant, the propagator explores about 47\,000 branches per second, against 2500 for the lazy variant. Table~\ref{tab:components} places these results in the context of the whole method. With lazy cuts, the CP phase already improves the mean upper bound from $+5.96\%$ to $-2.03\%$ relative to LBBD\_B\&B and proves 95 instances. The propagator further improves these results to $-3.62\%$ and 180 proven instances. The same conclusion holds on the pJSSP: the lazy variant proves 29 instances, fewer than Mistral (51), whereas \cpw{} proves 67. The advantage of \cpw{} therefore comes primarily from the combination of \texttt{PreemptiveNoOverlap} and the CP Optimizer search, rather than from the MILP bound or the starting solution. The lazy variant may exceed the time limit by a few seconds because it stops between consecutive searches (by at most 21\,s on the pJSSP), which, if anything, favors the lazy-cut approach.

\subsection{Starting solutions} \label{sec:res-start} Table~\ref{tab:start} evaluates the starting solutions on the 23 largest instances (Dauz\`ere-P\'er\`es and Hurink vdata la36--40), for which the impact of the initial solution is the greatest. We compare EDD on the MILP solution pool alone with the default strategy that also includes the window greedy of Section~\ref{sec:start}. \begin{table}[htbp] \centering\small \caption{Starting solutions on the 23 largest pFJSSP instances (600\,s). Start above final UB: mean gap of the starting solution to the final upper bound.} \label{tab:start} 
\resizebox{\textwidth}{!}
{
\begin{tabular}{l rrr r}
\toprule
Starting solution & UB vs LBBD (\%) & Start above final UB (\%) & Proven & Better / equal / worse$^a$ \\
\midrule
Pool only & $-$0.57 & 14.3 & 4 & -- \\
Pool + window greedy & $-$6.07 & 5.5 & 7 & 17 / 6 / 0 \\
\bottomrule
\end{tabular}
} 
\smallskip \begin{minipage}{0.95\textwidth}\footnotesize $^a$ Final upper bound compared with that of the previous line. \end{minipage} \end{table} The CP search is deterministic for a given starting solution, so the differences in Table~\ref{tab:start} are due solely to the initialization. The window greedy complements the MILP assignments, which balance machine loads but ignore sequencing, by choosing machines while scheduling. With the window greedy, the starting solution lies 5.5\% above the final upper bound on average, compared with 14.3\% for EDD alone. It improves the final upper bound on 17 instances and never worsens it. It also leads to three additional proofs of optimality, including v-la39 (922), previously bounded by $[922,949]$. On three instances, it already constructs an optimal schedule. Table~\ref{tab:heuristics} evaluates each starting heuristic separately on the full benchmark. To isolate the quality of the starts, only the first step of \cpw{} is executed: the MILP (with its 60\,s limit) followed by the starting heuristics, without the CP search. The gap is measured with respect to the best known upper bound, i.e., the best value among our runs, the published results, and the best known non-preemptive makespan. \begin{table}[htbp] \centering\small \caption{Starting heuristics, first step of \cpw{} alone: mean gap of the starting solution to the best known upper bound, number of instances on which the heuristic gives the best start (ties counted for each), and mean construction time. The EDD schedules use the assignments of the MILP; on the pJSSP, the assignment is fixed and the pool contains a single solution. The default start of \cpw{} is the best of the three.} \label{tab:heuristics} 
\begin{tabular}{l rrr rrr}
\toprule
 & \multicolumn{3}{c}{pFJSSP (276)} & \multicolumn{3}{c}{pJSSP (78)} \\
\cmidrule(lr){2-4}\cmidrule(lr){5-7}
Starting solution & gap (\%) & best & time (ms) & gap (\%) & best & time (ms) \\
\midrule
EDD, MILP incumbent & 19.27 & 96 & 0.02 & 16.65 & 68 & 0.03 \\
EDD, other pool assignments & 20.59 & 67/252 & 0.22 & -- & -- & -- \\
Window greedy (MWKR) & 14.16 & 145 & 7.18 & 20.06 & 16 & 13.75 \\
\midrule
Best of the three (default start) & 11.10 & -- & -- & 16.31 & -- & -- \\
\bottomrule
\end{tabular}
\end{table} 
No heuristic dominates the others. On the pFJSSP, the window greedy is the best individual heuristic, with a mean gap of 14.2\% to the best known upper bound and the best start on 145 instances. The EDD schedules built from the MILP incumbent and from the other pool assignments provide the best start on 96 and 67 instances, respectively (ties counted for each). On the pJSSP, where the machine assignment is fixed, the EDD schedule of the MILP solution gives the best start on 68 instances and the window greedy on 16. Taking the best solution produced by the three heuristics reduces the mean gap of the starting solution to 11.1\% on the pFJSSP and 16.3\% on the pJSSP. These values remain far from the final results of \cpw{}, confirming the conclusion of Table~\ref{tab:components}: starting solutions are useful, but most of the improvement comes from the CP search. All three heuristics require only a few milliseconds.

\section{Conclusions}
\label{sec:conclusion}
We have presented an exact fragment-free constraint programming framework for the preemptive flexible job shop scheduling problem (pFJSSP) with makespan minimization. By managing operation placement through flexible time windows rather than explicit processing segments, the method completely eliminates fragment variables and continuous-time discretization. Machine feasibility is enforced by extending the \texttt{PreemptiveNoOverlap} propagator to optional interval variables, performing overload checking and preemptive edge-finding in $O(n \log n)$ time per call while filtering infeasible machine alternatives. Initial lower bounds and starting solution pools are generated via a multi-round assignment relaxation incorporating Carlier's head--tail bounds, while Jackson's preemptive earliest-deadline-first rule reconstructs exact processing fragments for every incumbent solution in $O(n \log n)$ time.

Extensive computational evaluation across 276 pFJSSP and 78 pJSSP benchmark instances demonstrates that CPW substantially advances the state of the art. On the pFJSSP, CPW improves the best published upper bounds by 2.82\% on average, increases the number of proven optimal instances from 102 to 180 within 600~s, and closes 74 open instances for the first time, including the entire Barnes benchmark family. On the traditional pJSSP under identical hardware and time limits, CPW proves 67 out of 78 instances optimal, compared to 51 by the state-of-the-art solver Mistral, and closes 9 of the 11 long-standing open benchmark instances. Controlled component ablations further confirm the central role of the global constraint: replacing \texttt{PreemptiveNoOverlap} with a lazy overload cut separation scheme roughly halves the number of optimality proofs across both benchmarks (dropping from 180 to 95 on pFJSSP, and from 67 to 29 on pJSSP).

Three directions for future research emerge from our experiments. First, lower bounds can be strengthened by moving beyond single-machine relaxations. For the 72 open pFJSSP instances where the assignment relaxation solved to optimality within the time limit, tighter lower bounds will require energy-based overload reasoning over groups of candidate machines or destructive bounding executed directly inside the CP engine. For the 22 large instances where the MILP relaxation timed out with an average gap of 4.2\%, specialized Benders decomposition techniques for the assignment master problem offer a natural path forward. Second, on 94 open instances, the best known non-preemptive schedule remains lower than the best preemptive makespan found by CPW; integrating non-preemptive starting heuristics and targeted large-neighborhood search (LNS) operators could significantly enhance primal performance. Finally, the extension of \texttt{PreemptiveNoOverlap} to optional intervals provides a general mechanism that can be readily applied to other multi-resource preemptive domains, such as multi-mode resource-constrained project scheduling (MRCPSP).

\section*{Acknowledgements}
This research was supported by the Secretariat of Science, Humanities, Technology, and Innovation of Mexico through grant CBF-2025-I-2005.

\section*{Data availability statement}
The source code of \cpw{}, the benchmark instances, and the raw experimental results are available at: 
\url{https://github.com/VincentALBoyer/OR-Repository} in the folder Preemptive-Flexible-Job-Shop-Scheduling-Problem.

\bibliographystyle{plainnat}
\bibliography{references}

\end{document}